\documentclass[11pt]{amsart}
\usepackage{amssymb,latexsym}
\newtheorem{theorem}{Theorem}[section]
\newtheorem{prop}[theorem]{Proposition}
\newtheorem{lemma}[theorem]{Lemma}
\newtheorem{remark}[theorem]{Remark}
\newtheorem{question}[theorem]{Question}
\newtheorem{definition}[theorem]{Definition}
\newtheorem{cor}[theorem]{Corollary}

\newcommand{\C}{\mathbb{C}}

\begin{document}

\title{Symplectic birational geometry}

\author{ Tian-Jun Li$^1$ \& Yongbin Ruan$^2$}

\address{School  of Mathematics\\  University of Minnesota\\ Minneapolis, MN 55455}
\email{tjli@math.umn.edu}
\address{Department of Mathematics\\ University of Michigan\\ Ann Arbor, MI
48109-1109
 }\email{ruan@umich.edu}
 \thanks{${}^{1,2}$supported by NSF Grants}

\maketitle

\begin{center}
Dedicated to the occasion of Yasha Eliashberg's 60-th birthday
\end{center}

\section{Introduction}

Birational geometry has always been a fundamental topic in algebraic
geometry. In the 80's, an industry called Mori's birational geometry
program was created for the birational classification of algebraic
manifolds of dimension three.
   Roughly speaking, the idea of Mori's program
is to divide algebraic varieties into two categories: uniruled
versus non-uniruled. The  uniruled varieties are those containing a
rational curve through every point. Even for this class of algebraic
manifolds, the classification is usually not easy. So one is content
to carry out the classification of the more restrictive Fano
manifolds and prove some structural theorems such as the Mori fiber
space structure of any uniruled variety with Fano fibers. For
non-uniruled manifolds, one wishes to construct a ``minimal model"
by a sequence of contraction analogous to the blow-downs. One
immediate problem is that a contraction of smooth variety often
results a singular variety.  This technical problem often makes the
subject of birational geometry quite difficult. Only recently, the
 minimal model program was carried out to a large extent in
higher dimensions in the remarkable papers \cite{BCHM} and
\cite{Yiu}.

 In the early 90's, the second
author observed that some aspects of this extremely rich program of
Mori can be extended to symplectic geometry via the newly created
Gromov-Witten theory \cite{R1}. Specifically, he extended Mori's
notion of extremal rays to the symplectic category and used it to
study the symplectomorphism group. A few years later, Koll\'ar and
the second author showed that a smooth projective uniruled  manifold
carries a non-zero genus zero Gromov-Witten invariant with a point
insertion.  Shortly after, further deep relations between the
Gromov-Witten theory and birational geometry were discovered in
\cite{R2}, resulting in the speculation  that there should be a
symplectic birational geometry program.
    In the
mean time the Gromov-Witten theory, together with the Seiberg-Witten
theory, was applied with spectacular success to obtain basic
structure theorems of symplectic 4-manifolds, especially the
rational and ruled ones, cf. \cite{Mc1}, \cite{T}, \cite{LM}, \cite
{LL1}.

To be more specific, we define a genus 0 GW class as a nonzero
degree 2 homology class supporting nontrivial genus 0 GW invariants.
Let us then ask the following question: {\em what kind of structures
of a symplectic manifold are detected by its genus zero GW classes?}
What we would like to convey in this article is that the answer is
 precisely the {\em symplectic birational structure!}
 To start with consider  the sweep out of all the pseudo-holomorphic  rational curves in a GW class. The extreme case is
that the sweep out is the whole space. In this case the manifold is
likely a {\em uniruled manifold}. In general, the sweep out  is
likely a possibly singular {\em uniruled submanifold}. It has been
gradually realized that the symplectic birational geometry deals
precisely with these uniruled manifolds and uniruled submanifolds.

 In this article, we outline the main
elements of this new program in symplectic geometry. Let us first
mention that some technical
 difficulties in the algebraic birational geometry are still present
in our program, but we might be able to treat them with more
symplectic topological techniques. This is certainly the case for
contractions.
  Recall that the goal of birational geometry is to classify
algebraic
   varieties in the same birational class. Two
   algebraic varieties are birational to each other if and only if there is
    {\em a
   birational map} between them. A birational map is an isomorphism between  Zariski open
   sets, but it is not necessarily defined everywhere.
   If a birational map is defined everywhere, we
   call it a contraction.
   A contraction changes a lower dimensional uniruled subvariety only,
   hence we can view it as a surgery.
   Intuitively, a contraction simplifies a smooth variety,
   but as already mentioned, it often produces a singular variety. Various other types of surgeries are needed to deal with
   the resulting singularities. The famous ones are flip and flops which
   are much more subtle than contractions. We certainly can not avoid some aspects of this issue in our
program.
   A major problem in
   our program is then to see whether the flexibility in the symplectic
   category can produce many such kinds of surgery operations. In particular, we would like
    to interpret and construct
   flips and flops symplectically.

A new phenomenon in our program is that many obvious properties of
algebraic birational geometry are no longer obvious in the
symplectic category. Notably, the birational invariance of
uniruledness in \cite {HLR} is such an example, where we have to
draw newly developed powerful technology from the Gromov-Witten
theory.

However,  this perspective makes the subject distinctively
symplectic. And despite of these `old' and `new' obstacles, major
progress has been made recently in \cite{HLR}, \cite{tLR}, \cite
{tLR2}, \cite{Mc3}.

One eventual and remote goal of symplectic geometry is to
 classify symplectic manifolds.
 Symplectic birational geometry can   be considered as the first step towards such a
classification \cite{R2}. In addition, the author hope that
symplectic topological techniques and view points  in this new
program will also bring some fresh insight to the birational
classification of algebraic manifolds.

The article is organized as follows. We will set up symplectic
birational equivalence in section two. The transition as an extended
symplectic birational transformation will also be discussed. Section
three is devoted to the birational invariance of uniruled manifolds
and its classifications. In section four, we will discuss the
dichotomy of uniruled submanifolds. In section five, we will briefly
discuss speculations on minimal symplectic manifolds. We describe
various GW correspondences in section six. We finish the paper by
several concluding remarks. We should mention that this article does
not contain any proofs but provide appropriate references for the
results discussed.

Both author would like to thank Y. Eliashberg for the inspiration
and encouragement over so many years.  It is our pleasure to
dedicate this article to the occasion of his 60-th birthday. Special
thanks go to D. McDuff for many inspirational discussions and
careful readings of manuscripts of our papers as well as many useful
comments. Her influence on this subject is everywhere. We appreciate
the referee's suggestions which made the article more readable. We
are grateful to J. Hu for fruitful collaborations. Discussions with
Y. P. Lee and D. Maulik on recent GW techniques in algebraic
geometry, with J. Dorfmeister on relative symplectic cones, and with
Weiyi Zhang on the geometry of uniruled manifolds are really
helpful. We also appreciate the interest of V. Guillemin, H. Hofer,
C. LeBrun, M. Liu,
 N. Mok, C. Voisin, C. Taubes, S. T. Yau, A. Zinger and many others.
Finally we thank the algebraic geometry FRG group to organize a
wonderful conference on this subject in Stony Brook.

\section{Birational equivalence in symplectic geometry}

For a long time, it was not really clear what is an appropriate
notion of birational equivalence in symplectic geometry. Simple
birational operations such as blow-up/blow-down were known in
symplectic geometry for a long time \cite{GS, McS}. But there is no
straightforward generalization of  the notion of a general
birational map in the flexible symplectic category. The situation
changed a great deal when the weak factorization theorem was
established recently (see the lecture notes \cite{Ma} and the
reference therein) that any birational map between projective
manifolds can be decomposed as a sequence of blow-ups and
blow-downs. This fundamental result resonates perfectly with the
picture of the wall crossing of symplectic reductions analyzed by
Guillemin-Sternberg in the 80's. Therefore, we propose to use
their notion of cobordism in \cite{GS} as the symplectic analogue
of the birational equivalence (see Definition \ref{bc}). To avoid
confusion with other notions of cobordism in the symplectic
category, we would call it {\it symplectic birational cobordism}.

\subsection{Birational equivalence}

The basic reference for this section is \cite{GS}. We start with
the definition which is essentially contained in \cite{GS}.

\begin{definition} \label{bc} Two symplectic manifolds $(X, \omega)$ and
$(X', \omega')$ are birational cobordant if there are a finite
number of symplectic manifolds $(X_i, \omega_i)$,$ 0\leq i\leq k$,
with $(X_0, \omega_0)=(X, \omega)$ and $(X_k, \omega_k)=(X',
\omega')$, and for each $i$, $(X_i, \omega_i)$ and $(X_{i+1},
\omega_{i+1})$ are symplectic reductions of a semi-free Hamiltonian
$S^1$ symplectic manifold $W_i$ (of 2 dimensions higher).

\end{definition}

Here an $S^1$ action is called semi-free if it is free away from
the fixed point set.

There is  a related notion in dimension 4 in \cite{OO}. However we
remark that the cobordism relation studied in this paper is  quite
different from some other notions of symplectic cobordisms, see
\cite{EGH}, \cite{EH}, \cite{Gi}, \cite{GGK}.

According to \cite{GS}, we have the following basic factorization
result.

\begin{theorem} \label{f}
 A birational cobordism can be decomposed
as a sequence of elementary ones, which are modeled on blow-up,
blow-down and $\mathbb Z-$linear deformation of symplectic
structure.
\end{theorem}

A ${\mathbb Z}-$linear deformation is a path of symplectic form
$\omega+t\kappa$, $t\in I$, where $\kappa$ is a closed $2-$form
representing an integral class and $I$ is an interval. It was shown
in \cite{HLR} that ${\mathbb Z}-$linear deformations are essentially
the same as general deformations.

Observe that a polarization on a projective manifold, which is
simply a very ample line bundle, gives rise to a symplectic form
with integral class, well-defined up to isotopy. Together with the
weak factorization theorem mentioned in the previous page, we then
have

\begin{theorem}
Two birational projective manifolds with any polarizations are
birational cobordant as symplectic manifolds.
\end{theorem}

\section{Uniruled symplectic manifold}

\subsection{Basic definitions and properties}
Let us first recall the notion of uniruledness in algebraic
geometry.

\begin{definition} A projective manifold $X$ (over $\mathbb C$) is called (projectively) uniruled
if for every $x\in X$ there is a morphism $f:{\mathbb P}^1\to X$
satisfying $x\in f(\mathbb P^1)$, i.e. $X$ is covered by rational
curves.

\end{definition}

A beautiful property of of a uniruled projective manifold
(\cite{KMM}, \cite{K1}) is that general rational curves are
unobstructed, hence regular in the sense of the Gromov-Witten
theory. It is this property which underlies the aforementioned
result of Koll\'ar and Ruan (stated here in a sharper form noticed
by McDuff, cf. \cite{tLR}).

\begin{theorem} \label{projective2}A projective manifold is projectively uniruled if and only if
$$ \langle [pt], [\omega]^p, [\omega]^q\rangle^X_A> 0$$ for a nonzero class $A$, a K\"ahler form $\omega$
and integers $p, q$. Here $[pt]$ denotes the fundamental cohomology
class of $X$.

\end{theorem}

 In this article, for a closed symplectic manifold $(X, \omega)$, we denote its genus zero GW
invariant  in the curve class $A\in H_2(X;\mathbb Z)$ with
cohomology constraints $\alpha_1,\cdots, \alpha_k\in H^*(X;\mathbb
R)$ by \begin{equation}\label{inv}\langle \alpha_1, \cdots,
\alpha_k\rangle_A^X. \end{equation}  To define it we need to first
choose an $\omega-$tamed almost complex structure $J$ and consider
the moduli space of $J-$holomorphic rational curves with $k$ marked
points in the class $A$. Via the evaluation maps at the marked
points we pull back $\alpha_i$ to cohomology classes over the moduli
space, and the invariant (\ref{inv}) is supposed to be the integral
of the cup product of the pull-back classes over the moduli space.
Due to compactness and transversality issues the actual definition
requires a great deal of work (\cite{FO}, \cite{LT}, \cite{S1},
\cite{R4}). Intuitively the GW invariant (\ref{inv})  counts
$J-$holomorphic rational curves in the class $A$ passing through
cycles Poincar\'e dual to $\alpha_i$.

\begin{definition} Let $A\in H_2(X;{\mathbb Z})$ be a nonzero  class.   $A$ is called a GW class
if there is a non-trivial genus zero GW invariant of $(X, \omega)$
with curve class $A$. $A$ is said to be a uniruled class if it is a
GW class and moreover,  there is a nonzero GW invariant of the form
\begin{equation}\label{invariant}\langle [pt], \alpha_2, \cdots ,
\alpha_k\rangle^X_A,
\end{equation}
where $\alpha_i\in H^*(X;\mathbb R)$. $X$ is said to be
(symplectically) uniruled if there is a uniruled class.

\end{definition}

\begin{remark}
It is easy to see that  we could well use the  GW invariants with a
disconnected domain  to define this concept, subject to the
requirement that the curve component with the $[pt]$ constraint
represent nonzero class in $H_2(X;\mathbb Z)$. This flexibility is
important for the proof of the birational cobordism invariance.
\end{remark}

This notion has been studied in the symplectic context by G. Lu (see
\cite{Lu1}, \cite{Lu3}).  Notice that, by \cite{Ltj}, it is not
meaningful to define this notion by requiring that there is a
symplectic sphere in a fixed class through every point, otherwise
every simply connected manifold would be uniruled.

\begin{remark} \label{strong}According to Theorem \ref{projective2} a projectively
uniruled manifold is symplectically uniruled, in fact strongly
symplectically uniruled. Here $X$ is said to be strongly uniruled
if there is a nonzero invariant of the form (\ref{invariant}) with
$k=3$.
\end{remark}

Obviously the only uniruled 2--manifold is $S^2$.  In dimension 4
the converse is essentially true (see Theorem \ref{dim4}). While in
higher dimensions it follows from \cite{G} (see also \cite{Lu3})
that there are uniruled symplectic manifolds which are not
projective, and it follows from \cite{R3} that there could be
infinitely many distinct uniruled symplectic structures on a given
smooth manifold.

There are also  descendant GW invariants which are variations of the
GW invariants with constraints of the form $\tau_{j_i}(\alpha_i)$.
Here the class $\tau_{j_i}(\alpha_i)$ over the moduli space  is the
cup product of the pull-back of  the class $\alpha_i$ and the
$j_i-th$ power of a natural degree 2 class, which  is the 1st Chern
class of the (orbifold) line bundle over the moduli space whose
fibers are the cotangent lines at the $i-th$ marked point.

It is very useful to characterize uniruledness using these more
general  GW invariants (\cite{HLR}).

\begin{theorem} \label{descendant} A symplectic manifold $X$ is uniruled
if and only if there is a nonzero, possibly disconnected genus zero
descendant GW invariant
\begin{equation}\label{d-invariant}\langle \tau_{j_1}([pt]), \tau_{j_2}(\alpha_2), \cdots ,
\tau_{j_k}(\alpha_k)\rangle^X_A
\end{equation}
such that  the curve component with the $[pt]$ constraint  has
nonzero curve class.
\end{theorem}

 In particular, Theorem \ref{descendant}
is used to establish the fundamental birational invariance property
of uniruled manifolds in \cite{HLR}.

\begin{theorem}
Symplectic uniruledness is invariant under symplectic blow-up and
blow-down.
\end{theorem}


\subsection{Constructions}

An important aspect of symplectic birational geometry is the
classification of uniruled manifolds. This remains to be a distant
goal. A more immediate problem is to construct more examples.
Koll\'ar-Ruan's theorem shows that all the algebraic uniruled
manifold is symplectic uniruled. Another class of example is from
the following beautiful theorem of McDuff in \cite{Mc3}.
\begin{theorem}
Any Hamiltonian $S^1$-manifold is uniruled. Here a Hamiltonian
$S^1$-manifold is a symplectic manifold admitting a Hamiltonian
$S^1$-action.
\end{theorem}

A rich source of uniruled manifolds comes from
    almost complex uniruled  fibrations. Suppose that $\pi: X\rightarrow B$ is a fibration
    (with possibly singular fibers) where $X$ and $B$ are
    symplectic manifolds. We call it {\em an almost complex
    fibration} if there are tamed $J, J'$ for $X, B$ such that
    $\pi$ is almost complex.
    Symplectic fiber
    bundles over symplectic manifolds in the sense of Thurston are almost complex fibrations.
    Lefschetz fibrations, or more generally, locally holomorphic fibrations studied in \cite{Go2}
    are also almost complex.

    Let $\iota: \pi^{-1}(b)\rightarrow X$ be the embedding for a
    smooth  fiber over $b\in B$. We have the following result in
    \cite{tLR} by a direct geometric argument.

\begin{prop}\label{fiber}
    Suppose that $\pi: X\rightarrow B$ is an almost complex fibration between symplectic
    manifolds $X, B$. Then, for $A\in H_2(\pi^{-1}(b);\mathbb Z)$ and
    $\alpha_2,...,\alpha_k\in H^*(X;\mathbb R)$,
    \begin{equation}\label{=}
    <[pt], \iota^*\alpha_2, \cdots,
    \iota^*\alpha_k>^{\pi^{-1}(b)}_A=<[pt], \alpha_2, \cdots,
    \alpha_k>^X_{\iota_*(A)}.
    \end{equation}

    \end{prop}

    \begin{cor}\label{fibration}
Suppose that $\pi: X\rightarrow B$ is an almost complex fibration
between symplectic
    manifolds $X, B$. If a smooth fiber is uniruled and
    homologically
    injective (over $\mathbb R$), then $X$ is uniruled.
\end{cor}

The homologically injective assumption could be a strong one.
 Notice that for a fiber bundle,  the Leray-Hirsch
theorem asserts that, under the  homologically injective
assumption, the homology group of the total space is actually
isomorphic to the product of the homology group of the fiber and
the base. However, Corollary \ref{fibration} can still be applied
for all  product bundles, and
 all projective space fibrations (more generally, if the rational
cohomology ring of a smooth uniruled fiber is generated by
   the restriction of $[\omega]$).

Moreover, we were informed by McDuff that a Hamiltonian bundle is
homologically injective (or equivalently, cohomologically split)
if (cf. \cite{LM2})

a)  the base is $S^2$ (Lalonde-McDuff-Polterovich), and more
generally, a complex blow up of a product of projective spaces,

b) the fiber satisfies the hard Lefschetz condition (Blanchard),
or its real cohomology is generated by $H^2$.

 Here is another variation. As in the case of a projective space,
for a uniruled manifold up to dimension $4$,  insertions of a
uniruled class can all be assumed to be of the form $[\omega]^i$,
thus we also have

\begin{cor} If the general fibers of a possibly singular uniruled
fibration are $2$-dimensional or $4-$dimensional, then the total
space is uniruled.
\end{cor}

This in particular applies to a $2-$dimensional  symplectic conic
bundle. A symplectic conic bundle  is a  conic hypersurface bundle
in a smooth $\mathbb P^k$ bundle. Holomorphic conic bundles are
especially important in the theory of $3-$folds. It is conjectured
that a projective uniruled $3-$fold is either birational to a
trivial $\mathbb P^1-$bundle or a conic bundle.

Another important construction  first analyzed by McDuff is the {\em
divisor to ambient space} procedure.  It is part of the dichotomy of
uniruled divisors and would be discussed in  the next section (see
Theorem \ref{main1}).

\subsection{Geometry}\label{geometry}

Recall that the symplectic canonical class $K_{\omega}$ of $(X,
\omega)$ is defined to be $-c_1(TX, J)$ for any $\omega-$tamed
almost complex structure $J$.
 Observe that for a uniruled manifold $K_{\omega}$ is negative on any uniruled class by a simple
dimension computation of the moduli space. In particular,
$K_{\omega}$ cannot be represented by an embedded symplectic
submanifold. It leads to the following intriguing
question.\footnote{The converse of this question should be compared
with the Mumford conjecture: a projective manifold is uniruled if
and only if  it has Kodaira dimension $-\infty$.}

\begin{question} \label{sign} Does a uniruled manifold of (real) dimension $2n$ have a negative
$K_{\omega}^i\cdot [\omega]^{n-i}$ for some $i$?
\end{question}

On the other hand, we could ask if the canonical class $K_{\omega}$
is negative in the sense $K_{\omega}=\lambda [\omega]$ for
$\lambda>0$, is the symplectic manifold uniruled? Such a manifold is
known as a monotone manifold in the symplectic category, and it is
the analogue of a Fano manifold. Fano manifolds are projectively
uniruled by the famous bend-and-break argument of Mori. In fact,
Fano manifolds are even rationally connected.

We also observe here that uniruled manifolds satisfy a simple ball
packing constraint. To state it let us introduce the notion of a
minimal uniruled class, which is a uniruled class with minimal
symplectic area among all uniruled classes. This notion will also
play a crucial role in subsection \ref{non-negative}. Then it
follows from Gromov's monotonicity  argument  that the size of any
embedded symplectic ball is bounded by the area of a minimal
uniruled class.

\subsection{Rationally connected manifolds}
    A projective manifold is rationally connected if given any two
points $p$ and $ q$ there is a rational curve
 connecting $p$ and $ q$. It is equivalent to ``chain rational connected"
 where there is  a chain of rational curves connecting $p$ and $q$.

 The  outstanding conjecture is:
 {\em A projective manifold $X$ is rationally connected if and only if there
 is a nonzero {\it connected} GW invariant of the form
 \begin{equation} \label{rat-conn}<\tau_{j_1}([pt]), \tau_{j_2}([pt]), \cdots,
 \tau_{j_k}(\beta_k)>^X_{A}\neq 0\end{equation}
 for $A\neq 0$.}

We define a symplectic manifold to be rationally connected if there
is a nonzero invariant of the form (\ref{rat-conn}). Obviously, a
rationally connected manifold is uniruled.

Whether symplectic rational connectedness is a birational property
appears to be a hard question (cf. \cite{Vo}), though we believe it
is possible to show the invariance under certain types of blow-ups.
Characterizing such manifolds is likewise more difficult. But at
least we know that all such manifolds in dimension 4 are rational.
Moreover, it is expected  that symplectic manifold containing a
rationally connected symplectic divisor with certain strong
positivity is rationally connected (Initial progress has been made
in \cite{HRQ}).

Of course we can also similarly define $N-$rationally connectedness
for any integer $N\geq 2$. Moreover, it is not hard to see that
there is a parking constraint for $N$ balls in terms of (the area
of) a minimal $N-$rationally connected class.

\section{Dichotomy of uniruled divisors}

We have seen that up to birational cobordism symplectic manifolds
are naturally divided into uniruled ones and non-uniruled ones. In
this section we discuss uniruled submanifolds of codimension 2,
which we simply call symplectic divisors. One motivation comes from
the basic fact in algebraic geometry that
   various  birational surgery operations such as
   contraction and
    flop have a common feature: the subset being operated on
    is necessarily
    uniruled.

     Our key observation is that, as in the projective
    birational program, such a uniruled symplectic divisor
    admits a dichotomy  depending on the positivity of its normal
    bundle. If the normal bundle is non-negative in certain sense, it
    will force the ambient manifold to be uniruled.
    If the normal
    bundle is negative in certain sense, we can `contract' it to
    simplify the ambient  manifold.

    We have a rather general result in the non-negative case. In the
    second case our progress is limited to  simple uniruled divisors in
    6--manifolds.

\subsection{Dichotomy of uniruled divisor--non-negative
half}\label{non-negative}

     Suppose that $\iota: D\rightarrow (X, \omega)$
     is a  symplectic divisor (which we always assume to be smooth). Let $N_D$ be the normal bundle of
     $D$ in $X$. Notice that $N_D$ is a $2-$dimensional symplectic vector bundle and
     hence has a well defined first Chern class.  We will often  use $N_D$ to denote the first Chern class.

\begin{theorem} \label{main1} Suppose $D$ is uniruled and $A$ is a minimal uniruled class of $D$
such that
\begin{equation}\label{k}
<\iota^*\alpha_1, \cdots, \iota^*\alpha_l, [pt], \beta_2, \cdots,
\beta_k>^D_{A}\neq 0
\end{equation} for $\alpha_i\in H^*(X;\mathbb R), \beta_j\in H^*(D;\mathbb R)$ with $\beta_1=[pt]$, and $k\leq N_D(A)+1$. Then $(X,\omega)$ is uniruled.
\end{theorem}

Here what matters in (\ref{k}) is the number of insertions which do
not come from $X$. There are situations where can simply take $k=1$
hence only require $N_D(A)$ be non-negative. In particular, we have

\begin{cor} \label{homologically injective} Suppose $D$ is a homologically injective uniruled
divisor of $X$ and the normal bundle $N_D$ is non-negative on a
minimal uniruled class. Then $X$ is uniruled.

\end{cor}

We can ask whether  the  converse of Theorem \ref{main1} is also
true. It is obvious
 in dimension 4. In higher dimension we could easily construct
 non-negative singular uniruled divisors. The hard and important
 question is whether we can they can be smoothed inside $X$.



As  mentioned in the previous section, Theorem \ref{main1} can also
be
    viewed as a construction of uniruled manifolds, generalizing
    several early results of McDuff. We list some examples here,
    more examples can be found in \cite{tLR}.

{\bf 4--dimensional uniruled divisors}: A deep result in dimension 4
is  that uniruled manifolds can be completely classified.

\begin{theorem} \label{dim4}  (\cite{Mc}, \cite{LL1}, \cite{LL2}, \cite{LM}, \cite{T}) A $4-$manifold $(M, \omega)$ is uniruled if and only if
it is rational or ruled.  Moreover, the isotopy class of $\omega$ is
determined by $[\omega]$.
\end{theorem}

 Here  symplectic $4-$manifold $(M, \omega)$ is called rational if its underlying smooth manifold $M$
 is $S^2\times S^2$ or
$\mathbb P^2\# k{\overline{\mathbb P}^2}$ for some non-negative
integer $k$. $(M, \omega)$  is called ruled if its underlying smooth
manifold $M$ is the connected sum of a number of (possibly zero)
${\overline {\mathbb P}^2}$ with an $S^2-$bundle over a Riemann
surface.

We need to analyze minimal uniruled classes and the corresponding
insertions.

\begin{prop} \label{criterion} If $A$ is a uniruled class of a $4-$manifold, then
$A$ is represented by an embedded symplectic surface, and $A$
satisfies (i) $K_{\omega}(A)\leq -2$, (ii) $A^2\geq 0$, (iii)
$A\cdot B\geq 0$ for any class $B$ with a non-trivial GW invariant
of any genus.

\end{prop}

 For  $\mathbb P^2$, let $H$ be the generator of $H_2$ with
 positive area. $H$ is a uniruled class and any uniruled class
 of the form $aH$ with $a>0$.
 Obviously,  $H$ is the minimal uniruled class. The
relevant insertion is $([pt], [pt])$. As $[pt]$ is a restriction
class, i.e. an $\alpha$ class, we can take $k=1$.

Similarly, for the blow-up of an $S^2-$bundle over
 a surface of positive genus, the fiber class is a uniruled class,
 and any uniruled class is a positive multiple of the fiber class.
The relevant insertion for the fiber class  is $[pt]$. Thus again we
can take $k=1$.

It is easier to apply Theorem \ref{main1} in this case.

\begin{cor} Suppose $(X^6, \omega)$ contains a divisor $D$
 which is diffeomorphic to $\mathbb P^2$ or the blow-up of an $S^2-$bundle over
 a surface of positive genus. If the normal bundle $N_{D}$ is
 non-negative on
 a uniruled class, then $X$ is uniruled.
\end{cor}

 For other $M^4$, the uniruled classes are
not proportional to each other. Thus the minimality condition
depends on the class of the symplectic form on  $M$.

 We
first analyze  the easier case of an $S^2-$bundles over $S^2$. For
$S^2\times S^2$, by uniqueness of symplectic structures,  any
symplectic form is of product form. Let $A_1$ and $A_2$ be the
classes of the factors with positive area. It is easy to see that
any uniruled class is of the form $a_1A_1+b_1A_2$ with $a_1\geq 0,
a_2\geq 0$. Thus either $A_1$ or $A_2$ has the minimal area.

For the nontrivial bundle $S^2\tilde \times S^2=\mathbb P^2 \#
\overline {\mathbb P}^2$, let $F_0$ be the class of a fiber with
positive area and $E$ be the unique $-1$ section class with positive
 area. If $aF_0+bE$ is a uniruled class then $b\geq 0$ by (iii) of
Proposition \ref{criterion}, since $F_0\cdot E=1, F_0\cdot F_0=0$.
And if $b> 0$, then $a\geq 1$ by (i) of Proposition \ref{criterion}.
Thus $F_0$ is always the minimal uniruled class no matter what the
symplectic structure is.

Since the relevant insertion  for $A_1, A_2$ and $F_0$ is  just
$[pt]$, we have

\begin{cor}
Suppose $D=S^2\times S^2$ and the  restriction of the normal bundle
$N_{D}$  to the factor with the least area  is non-negative, then
$X$ is uniruled. In the case of the non-trivial bundle,
 $X$ is uniruled if the
restriction of the normal bundle  $N_{D}$ to $F_0$ is
non-negative.
\end{cor}

The remaining $M^4$ are connected sums of $\mathbb P^2$ with at
least 2 ${\overline {\mathbb P}^2}$. It is complicated to analyze
minimal uniruled classes in this case. In \cite{tLR} it is shown
that they are generated by the so called fiber classes.


{\bf Higher dimensional case: }
     In higher dimension, it is still a remote goal to classify
     all the uniruled symplectic manifolds. Instead of considering
     an arbitrary uniruled symplectic divisor, we start with
Fano hypersurfaces.

     When the divisor $D\subset (X, \omega)$ is symplectomorphic to a divisor of $\mathbb P^n$ (for $n\geq 4$) of degree at most
     $n$, $D$ is Fano and hence uniruled.
      Of course a particular case is $D=\mathbb P^{n-1}$. Since $n\geq 4$, by the Lefschetz
     hyperplane theorem, $b_2=1$ for $D$.
     According to Theorem \ref{projective2}, for a minimal
uniruled class $A$, we can take $k$ to be equal to $1$. Hence
 $X$ is uniruled if
       $N_D=\lambda [\omega|_D]$ with $\lambda\geq 0$.

In general case we still need to verify the minimal condition. Of
course the uniruled divisor needs not to be a projective manifold.
For instances, the divisor could be a rather general uniruled
fibration discussed in \S2. Let us treat the case of a symplectic
$\mathbb P^{k}-$bundle. Since the line class in the fiber is
uniruled, and the relevant insertions can be taken to be $([pt],
[\omega|_D]^k)$, we have

\begin{cor} Suppose $D$ is a symplectic divisor of $X$. If $D$ is a projective space bundle with
the fiber class being  the minimal uniruled class and  normal bundle
$N_D$  being non-negative
 along the fibers, then $X$ is unruled.
\end{cor}

McDuff also considered the case of product $\mathbb P^k-$bundles in
\cite{Mc2}. A natural source of such a $D$ is from blowing up a
`non-negative' $\mathbb P^k$ with a large trivial neighborhood.
Suppose $\mathbb P^k\subset X$
 has trivial normal bundle. Then the blow up along $\mathbb P^k$
has a divisor $D=\mathbb P^k\times \mathbb P^{n-k-1}$. The normal
bundle of $D$ along a line in $\mathbb P^k$ is trivial. Similar to
the case of $S^2\times S^2$, we can argue that
 the area of this line is minimal among all uniruled class
of $D$. In particular, as observed by \cite{Mc2},  a symplectic
$\mathbb P^1$ with a sufficiently large product symplectic
neighborhood can only exist in a uniruled manifold.

In fact we can prove more.

\begin{cor} \label{size} Suppose $S$ is a uniruled symplectic submanifold
whose minimal uniruled class has area $\eta$ and insertions all
being  restriction classes. If  $S$ has a trivial symplectic
neighborhood of radius at least $\eta$. Then $X$ is uniruled.
\end{cor}

\subsection{Symplectic `blowing-down' in dimension
six}\label{contraction}  Blowing up in dimension 6 gives rise to a
symplectic $\mathbb P^2$ with normal $c_1=-1$ or a symplectic
$S^2-$bundle with normal $c_1=-1$ along the $S^2-$fibers. A natural
question is whether such a uniruled divisor always arises from a
symplectic blow-up. In other words, we are interested in a criterion
for blowing-down. A nice feature is that by Theorem \ref{dim4} the
answer would also only depend on $[\omega]$.
 Moreover, as every symplectic structure on such a
4--manifold is K\"ahler, we can apply algebro-geometric techniques
to understand this problem.


The case of $\mathbb P^2$ is simple: it can always be blown down
just as in the case of $\mathbb P^1$ with self-intersection $-1$ in
a symplectic 4--manifold. For the case of an $S^2-$bundle over a
Riemann surface $\Sigma_g$ of genus $g$, it is more complicated and
perhaps more interesting. Topologically blowing down an $S^2-$bundle
over $\Sigma_g$ is the same as topologically fiber summing  with the
pair of a $\mathbb P^2-$bundle and an embedded $\mathbb P^1-$bundle
over $\Sigma_g$ with opposite normal bundle (see \ref{cut}). To
perform the fiber sum symplectically we also need to match the
symplectic classes of the divisors.
  For this purpose we need to understand the
relative symplectic cone of such a pair. We determine in \cite{tLR2}
the relative K\"ahler cone for various complex structures coming
from stable and unstable rank 3 holomorphic bundles over $\Sigma_g$.
Consequently we obtain


\begin{theorem}\label{fibersum} Let $(X, \omega)$ be a symplectic manifold of
dimension 6. Let $D$ be a symplectic divisor which is an
$S^2-$bundle over $\Sigma_g$. Suppose $N_D(f)=-1$ and $N_D(s)=d$,
where $f$ is the fiber class of $D$ , and $s$ is a section class of
$D$  with square $0$ if $D$ is a trivial bundle and square $-1$ if
$D$ is a non-trivial bundle. Further assume that $[\omega|_D](f)=a$
and $[\omega|_D](s)=b$. Then a symplectic fiber sum can be performed
with a symplectic $(\mathbb P^2, \mathbb P^1)/\Sigma_g$ pair if
either $d\geq 0$, or

(i) $d<0, g=0, b>[\frac{-d+2}{3}]a$,

(ii) $d<0, g\geq 1, b>\frac{-d}{3}a$.

\end{theorem}

Here $[x]$ denotes the largest integer bounded by $x$ from above.
This result is optimal in the case of genus $0$. For instance, we
show that a symplectic $S^2\times S^2$ with normal $c_1=-1$ along
each family of $S^2$ can be blown down if the symplectic areas of
the two factors are not the same. This picture is consistent with
the flop operation for projective 3--folds. On the other hand, when
$g\geq 1$, a symplectic $S^2\times \Sigma_g$ with normal $c_1=-1$
along each factor can be fiber summed if the symplectic area of the
$S^2$ factor is at most 3 times of that of the $\Sigma$ factor. We
would very much like to find out whether the restriction on the
areas is really necessary in the case of positive genus. The picture
would be rather nice if the restriction can eventually be removed.
On the other hand, it would be surprising, even intriguing, if it
turns out there is an obstruction.
 We are
also interested to see how much of the 6--dimensional investigation
can be carried out to higher dimensions.

Another remaining issue is whether such a symplectic fiber sum is
actually a birational cobordism operation. This is because that a
symplectic blowing down further requires the symplectic $(\mathbb
P^2, \mathbb P^1)/\Sigma_g$ pair have an infinity symplectic
section, which is a symplectic surface of genus $g$. We are
investigating wether our more general fiber sums are equivalent to
symplectic blowing down up to deformation. Either positive or
negative answer would be very interesting.

\section{Minimal symplectic manifolds}

\subsection{Minimality} Motivated by the Mori program for algebraic
3--folds and understandings of symplectic 4--manifolds (and
2--manifolds), we discuss in this section the notion of minimal
manifolds in dimension 6.

Let us first recall the notion of minimality by McDuff in dimension
4. Let ${\mathcal E}_X$ be the set of homology classes which have
square $-1$ and are represented by smoothly embedded spheres. We say
that $X$ is smoothly minimal if ${\mathcal E}_X$ is empty. Let
${\mathcal E}_{X,\omega}$ be the subset of ${\mathcal E}_X$ which
are represented by embedded $\omega-$symplectic spheres. We say that
$(X,\omega)$ is
 symplectically minimal if ${\mathcal E}_{X,\omega}$ is empty.
When $(X,\omega)$ is non-minimal,  one can blow down some of the
symplectic $-1$ spheres to  obtain a minimal symplectic $4-$manifold
$(N, \mu)$, which is called a (symplectic) minimal model of
$(X,\omega)$ ([Mc]). We summarize the basic facts about the minimal
models in the following proposition.

\begin{theorem} [\cite{Ltj4}, \cite{LL2}, \cite{Mc}, \cite{T}]\label{min}
Let $X$ be a  closed oriented smooth $4-$manifold and $\omega$ a
symplectic form on $X$ compatible with the orientation of $X$.

\noindent 1. $X$ is smoothly minimal if and only if $(X,\omega)$ is
symplectically minimal. In particular the underlying smooth manifold
of the (symplectic) minimal model of $(X,\omega)$ is smoothly
minimal.

\noindent 2.  If $(X,\omega)$ is not rational nor ruled, then it has
a unique (symplectic) minimal model. Furthermore, for any other
symplectic form $\omega'$ on $X$ compatible with the orientation of
$X$, the (symplectic) minimal models of $(X,\omega)$ and
$(X,\omega')$ are diffeomorphic as oriented manifolds.

\noindent 3. If $(X,\omega)$ is rational or ruled, then its
(symplectic) minimal models are diffeomorphic to $CP^2$ or an
$S^2-$bundle over a Riemann surface.
\end{theorem}

\begin{definition} \label{min6} A symplectic 6--manifold is {\em
       minimal} if it does not contain any  rigid stable uniruled
       divisor.
       \end{definition}
 Here, a uniruled divisor is stable if one
       of its uniruled classes $A$ is a GW class of the ambient manifold with $K_{\omega}(A)\leq -1$.
       And a uniruled divisor is rigid if none of its uniruled class
       is a uniruled class of the ambient manifold.

We observe this definition also applies to manifolds of dimensions 2
and 4. Every 2--manifold is obviously minimal as the only divisors
are points. For a 4 manifold any uniruled divisor is $S^2$. Let $A$
be the class of  a stable $S^2$. Then $K_{\omega}(A)\leq -1$. On the
other hand a rigid $S^2$ must have $K_{\omega}(A)<0$ by \cite{Mc}.
The only rigid stable uniruled divisor is an $S^2$ with
$K_{\omega}(A)=-1$. By the adjunction formula it is a symplectic
$-1$ sphere. Thus this definition of minimality  agrees with the
usual notion of minimality by McDuff.

An immediate question is whether any  6--manifold  has a minimal
model. A general strategy is to construct a minimal model by
performing consecutive
       contractions as discussed in \ref{contraction}.
  Notice that we do not need to eliminate  all negative uniruled submanifolds, only
those which are rigid and stable. Thus we might be able to avoid
complicated singularities. However, we will definitely encounter
orbifold singularities. The first author showed in \cite{R1} that
all ($K-$negative) extremal rays for algebraic 3--folds give rise to
non-trivial GW classes. One such an extremal ray arises from the
line class of a  $\mathbb P^2$ divisor with normal $c_1=-2$. To
carry out the contraction one has to enlarge to the category of
symplectic orbifolds.


Now we single out an important class of minimal manifolds.

\begin{definition}
We define a cohomology class $\alpha\in H^2(X,{\mathbb Z})$ to be
nef if it is non-negative on all GW classes.
\end{definition}

\begin{lemma} Let $(X, \omega)$ be a manifold with nef $K_{\omega}$. Then
$(M, \omega)$ is non-uniruled and minimal.

\end{lemma}
\begin{proof} The first statement is obvious as any uniruled class
$A$ of $X$ satisfies $K_{\omega}(A)\leq -2$.

There can not be any stable uniruled divisor in $X$ as
$K_{\omega}(A)\leq -1$ for class $A$ which is a uniruled class of a
stable divisor as well as a GW class of  $X$.

\end{proof}

A minimal model of a uniruled manifold is still uniruled and hence
can not be $K_{\omega}-$nef. The natural question is whether any
minimal model of a non-uniruled manifold must have  nef
$K_{\omega}$. The first step towards this question is to show that
any GW class $A$ with $K_{\omega}(A)\leq -1$ is a uniruled class of
a (smooth) divisor. Here the issue is again smoothing  singular
uniruled divisors.

A further  question is whether birational cobordant $K_{\omega}-$nef
manifolds are related by an analogue of the
       $K-$equivalence. Two algebraic manifolds $X, X'$ are
       $K-$equivalent if there is a common resolutions $\pi_1:
       Z\rightarrow X, \pi_2: Z\rightarrow X'$ such that $\pi^*_1
       K_X=\pi^*_2 K_{X'}$. $K-$equivalent manifolds have many
       beautiful properties, in particular,  they have the same betti
       numbers.


\subsection{Kodaira dimension}
     The notion of   Kodaira dimension has been defined for symplectic manifolds up to
       dimension 4 (\cite{Ltj2}, \cite{McS}). Whenever it is defined it is a finer invariant of
birational cobordism then uniruledness.

The Kodaira dimension of  a 2--dimensional symplectic manifold $(F,
\omega)$ is
 defined as
\[ \kappa(F,\omega)=\left\{ \begin{array}{ll}
-\infty &\hbox{if $K_{\omega}\cdot [\omega]<0$,} \\
0&\hbox{if $K_{\omega}\cdot [\omega]=0$,}\\
1&\hbox{if $K_{\omega}\cdot [\omega]> 0$.}
\end{array}
\right. \] Clearly $\kappa(F, \omega)=-\infty, 0, 1$ if and only if
the genus of $F$ is $0, 1, \geq 2$ respectively.  Notice that
$\kappa(F, \omega)=-\infty$ if and only if $(F, \omega)$ is
uniruled.

 For a minimal symplectic $4-$manifold $(X, \omega)$ its Kodaira dimension
 is defined in the following way (\cite{Le2}, \cite{McS}, \cite{Ltj2}):
\[ \kappa(X,\omega)=\left\{ \begin{array}{ll}
-\infty &\hbox{if $K_{\omega}\cdot [\omega]<0$ or $K_{\omega}\cdot K_{\omega}<0$},\\
0&\hbox{if $K_{\omega}\cdot [\omega]=0$ and $K_{\omega}\cdot K_{\omega}=0$},\\
1&\hbox{if $K_{\omega}\cdot [\omega]> 0$ and $K_{\omega}\cdot K_{\omega}=0$},\\
2&\hbox{if $K_{\omega}\cdot [\omega]>0$ and $K_{\omega}\cdot K_{\omega}>0$}.\\
\end{array}
\right. \]
 The Kodaira dimension of a non-minimal manifold is
defined to be that of any of its minimal models.

Based on the  Seiberg-Witten theory and properties of minimal models
(cf. Theorem \ref{min}, \cite{LL2}, \cite{Liu1}, \cite{OO}),
\cite{T}) it is shown in  \cite{Ltj2} that the Kodaira dimension
$\kappa(M,\omega)$ is well defined. In particular, we need to check
that a minimal 4-manifold cannot have
\begin{equation} \label{well defined} K_{\omega}\cdot
[\omega]=0, K_{\omega}\cdot K_{\omega}>0.
\end{equation}

We list some basic properties of $\kappa(X, \omega)$.
 It is also observed in \cite{Ltj2} that, if
$\omega$ is a K\"ahler form on a complex surface $(X, J)$, then
$\kappa(X, \omega)$ agrees with the usual holomorphic Kodaira
dimension of $(X, J)$.

$(X, \omega)$ has $\kappa=-\infty$ if and only if it is uniruled.

It is further shown in \cite{Ltj2}  that minimal symplectic
$4-$manifolds with $\kappa=0$ are exactly those with torsion
canonical class, thus they can be viewed as {\it symplectic
Calabi-Yau surfaces}. Known examples of symplectic $4-$manifolds
with torsion canonical class are either K\"ahler surfaces with
(holomorphic) Kodaira dimension zero or $T^2-$bundles over $T^2$. It
is shown in \cite{Ltj3} and \cite{B} that  a {\it minimal}
symplectic $4-$manifold with $\kappa=0$ has the rational homology as
that of K3 surface, Enriques surface or a $T^2-$bundle over $T^2$.

Suppose $(X, \omega)$ is a minimal 6--dimensional manifold.
 we propose to define its Kodaira dimension
 in the following way\footnote{Compare with Question \ref{sign}}:

\[ \kappa(X,\omega)=\left\{ \begin{array}{ll}
-\infty &\hbox
 {if one of
 $K_{\omega}^i\cdot [\omega]^{3-i}$ is negative,} \\
k &\hbox{if $K_{\omega}^i\cdot [\omega]^{3-i}=0$ for $i\geq
 k$, and $K_{\omega}^i\cdot [\omega]^{3-i}>0$ for $i\leq k$.}\\
 \end{array}
\right. \]

 Notice that, as in dimension 4, there is the issue of  well definedness of $\kappa(X, \omega)$.
 And this leads to some possible intriguing properties of
of minimal 6--manifolds, one of which is
 whether there is any minimal 6--manifold
with
 \begin{equation} \label{well defined} K_{\omega}\cdot
[\omega]^2=0, K_{\omega}^2[\omega]=0, K_{\omega}^3>0. \end{equation}

       \section{Correspondences in the Gromov-Witten theory}

       As the Gromov-Witten theory is built into the foundation of symplectic birational geometry,
       it is natural that we use many techniques from the Gromov-Witten theory such
       as localization and degenerations. It turns out that we have to
       use very sophisticated Gromov-Witten machinery. Take the birational invariance of the uniruledness as an example.
       The definition of uniruledness
       requires only a single non-vanishing GW invariant. However, it is well-known
       that a single Gromov-Witten invariant tends to  transform in a rather complicated
       fashion. On the other hand, it is often easier to control the transformation of the
        Gromov-Witten theory as a whole.
       One often phrases such a amazing phenomenon as a kind of correspondences.
       There are many examples such as the Donaldson-Thomas/Gromov-Witten correspondence \cite{MMOP}, the crepant
        resolution conjecture \cite{R3} and so on.

       The correspondences appearing in our context are not as strong as the above ones.
       Its first example is  the ``relative/absolute correspondence"
      constructed by
     Maulik-Okounkov-Pandharipande (\cite{MP}). It is the  generalization to the
     situation of blow-up/down by the authors and Hu which underlies the birational invariance
     of uniruledness. And Theorem \ref{main1} is proved by another
     technical variation of the relative/absolute correspondences incorporating divisor invariants.
     Roughly speaking, a correspondence in this context is a
     package to organize the degeneration formula in a very nice
     way.

We restrict ourselves to genus GW invariants in this article.

     \subsection{Symplectic
cut and the degeneration formula}\label{cut}

\subsubsection{Symplectic cut along a submanifold}
Let $(X,\omega)$ be a closed symplectic manifold. Let $S$ be a
hypersurface having a neighborhood with a free Hamiltonian
$S^1-$action. For instance, if there is a symplectic submanifold in
$X$, then the hypersurfaces corresponding to sphere bundles of the
normal bundle have this property. Let $Z$ be the symplectic
reduction at the level $S$, then $Z$ is the $S^1-$quotient of $S$
and is a symplectic manifold of 2 dimension less.

We can cut $X$ along $S$ to obtain two closed symplectic manifolds
$(\overline X^+,\omega^+)$ and $(\overline X^-,\omega^-)$ each
containing a smooth copy of $Z$, and satisfying $\omega^+\mid_Z =
\omega^-\mid_Z$ (\cite{Le}).

 In particular, the pair $(\omega^+,
\omega^-)$ defines a cohomology class of
$\overline{X}^+\cup_Z\overline{X}^-$, denoted by
$[\omega^+\cup_Z\omega^-]$.
 Let $p$ be the continuous collapsing map
$$p:X\to \overline{X}^+\cup_Z\overline{X}^-.$$
It is easy to observe that
\begin{equation}\label{cohomology relation}
   p^* ([\omega^+\cup_Z\omega^-]) = [\omega].
\end{equation}

Let $\iota:D\to X$ be a smooth connected symplectic divisor. Then we
can cut along $D$, or precisely, cut along a small circle bundle $S$
over $D$ inside $X$.

In this case, as a smooth manifold, $\overline{X}^+ =X$, which we
will denote by $\tilde X$. Denote the symplectic reduction of $S$ in
$\tilde X$ still by $D$. Notice however, the symplectic structure is
different from the original divisor. And $\overline{X}^-= \mathbb
P(N_D\oplus \underline {\mathbb C})$, the projectivization of
$\mathbb P(N_D\oplus \underline {\mathbb C})$\footnote{Notice that
our convention here is opposite to that in \cite{HLR}}. We will
often denote it simply by $P_D$ or $P$. Notice that $\mathbb
P(N_D\oplus \underline {\mathbb C})$ has two natural sections,
$$  D_{0}=\mathbb P(0\oplus \underline \C), \quad D_{\infty}=\mathbb P(N_D\oplus 0).$$
 The symplectic reduction of $S$ in $P_D$ is the
section $D_{\infty}$.

 In
summary, in this case,
 $X$ degenerates into $(\tilde X, D)$ and
$(P_D, D_{\infty})$. We also denote $\omega^-$ by $\omega_P$.

More generally, we can cut along a symplectic submanifold $Q$ of
codimension $2k$, or precisely, cut along a sphere bundle $S$ over
$Q$. Then  $\overline{X}^+$ is a symplectic blow up of $X$ along $Q$
and the symplectic reduction $Z\subset \overline{X}^+$ is the
exceptional symplectic divisor. In this case $\overline{X}^-=
\mathbb P(N_Q\oplus \underline {\mathbb C})$, which is a $\mathbb
P^k$ bundle over $Y$.

\subsubsection{Degeneration formula}\label{df}
 Given a symplectic cut, there is a  basic link between absolute invariants
 of $X$ and relative invariants
of $(\overline{X}^{\pm}, Z)$ in \cite{LR} (see also \cite{IP}, and
\cite{Li2} in  algebraic geometry).  We now  describe such a
formula.

 Let $B\in
H_2(X;{\mathbb Z})$ be in the kernel of
$$
  p_* : H_2(X;{\mathbb Z})
\longrightarrow H_2(\overline{X}^+\cup_Z\overline{X}^-; {\mathbb
Z}). $$
 By (\ref{cohomology relation}) we have $\omega(B) =0$.
Such a class is called a vanishing cycle.
 For $A\in H_2(X; {\mathbb Z})$ define $[A] = A + \mbox{Ker}
(p_*)$ and
\begin{equation}\label{vanishing cycle}
\langle\tau_{d_1}\alpha_1, \cdots, \tau_{d_k}\alpha_k\rangle^X_{[A]}
= \sum_{B\in[A]}\langle\tau_{d_1}\alpha_1, \cdots,
\tau_{d_k}\alpha_k\rangle^X_{B}.
\end{equation}

At this stage  we need to assume that each cohomology class
$\alpha_i$ is of the form \begin{equation}\label{distribution}
\alpha_i = p^*(\alpha_i^+\cup_Z\alpha_i^-).
\end{equation}
 Here $\alpha_i^\pm \in H^*(\overline{X}^\pm; {\mathbb
R})$ are classes with  $\alpha_i^+\mid_Z = \alpha_i^-\mid_Z$ so that
 they give rise to a class $\alpha_i^+\cup_Z\alpha_i^-\in
H^*(\overline{X}^+\cup_Z\overline{X}^-; {\mathbb R})$.

The degeneration formula expresses $\langle\tau_{d_1}\alpha_1,
\cdots, \tau_{d_k}\alpha_k\rangle^X_{[A]} $ as a sum of products of
relative invariants of $(\overline{X}^+, Z)$ and $(\overline{X}^-,
Z)$, possibly with disconnected domains.  In each product of
relative invariants, what is relevant for us are the following
conditions:

$\bullet$ the union of two domains  along relative marked points is
a stable genus $0$ curve with $k$ marked points,

$\bullet$ the total curve class  is equal to $p_*(A)$,

$\bullet$ the relative insertions are dual to each other,

$\bullet$ if $\alpha_i^+$ appears for $i$ in a subset of
$\{1,\cdots, k\}$, then $\alpha_j^-$ appears for $j$ in the
complementary subset of $\{1,\cdots, k\}$.

In the case of cutting along a symplectic submanifold it is easy to
show that all the invariants on the right hand side of
(\ref{vanishing cycle}) vanish except $$\langle\tau_{d_1}\alpha_1,
\cdots, \tau_{d_k}\alpha_k\rangle^X_{A}.$$ Thus the degeneration
formula computes $\langle\tau_{d_1}\alpha_1, \cdots,
\tau_{d_k}\alpha_k\rangle^X_{A} $ in terms of relative invariants of
$(\overline X^{\pm}, Z)$.

\subsection{Absolute/Relative, blow-up/down
and divisor/ambient space correspondences} Giving a symplectic
manifold $(X, \omega)$ we are interested in determining its GW
classes and uniruled classes. Suppose $X$ has some explicit
symplectic submanifolds, then we could cut $X$ and attempt to apply
the degeneration formula to compute a given GW invariant. However,
this is often impractical, as we need to know all the relevant
relative invariants of $(\overline X^{\pm}, Z)$ and relative
invariants are generally harder to compute themselves.

Remarkably it is shown in \cite{MP} that, in case the submanifold
$D$ is a divisor and hence $\overline X^+=X$, the degeneration
formula can be inverted to express a (non-descendant) relative
invariant of $(X, D)$ in terms of invariants of $X$ and relative
invariants of $(P_D, D_{\infty})$. We brief describe the strategy of
proof in \cite{MP}.

The first idea is to associate a possibly descendant  invariant of
$X$ to each non-descendant relative invariant of $(X, D)$ where
absolute insertions are kept intact and  contact orders of relative
insertions are replaced by appropriate descendant powers.

Observe then relative GW invariants are linear on the insertions. So
we can choose a generating set $I$ of non-descendant relative GW
invariants by choose bases of cohomology of $X$ and $D$ and require
the absolute and relative insertions lie in the two bases.

 The next idea is to introduce a partial order
on $I$ with 2 properties. Firstly, given a relative invariant of
$(X, D)$, when applying the degeneration formula to the associated
invariant of $X$, the given relative invariant is the largest one
among those relative invariants of $(X, D)$ appearing in the formula
and with nonzero coefficient. Recall that the right hand side of the
degeneration formula is a sum of products, for each product the
relative invariant of $(P_D, D_{\infty})$ is considered to be the
coefficient. Secondly, the partial order is lower bounded in the
sense there are only finitely many invariants in $I$ lower than any
given relative invariant in $I$.

Then inductively, any relative invariant in $I$ can be expressed in
terms of invariants of $X$ and relative invariants of $(P_D,
D_{\infty})$.

In \cite{HLR} we slightly reformulate the absolute/relative
correspondence as a lower bounded and triangular (and hence
invertible) transformation $T$ in an infinity dimension vector
space, sending the relative vector $v_{rel}^I$ determined by all
relative invariants in $I$ to the $v_{abs}^I$ absolute vector
determined by all the associated invariants of $X$. Furthermore, if
$I_{pt}$ is the subset of $I$ such that one of the absolute
insertions is a $[pt]$ insertion, $T$ still interchanges the
absolute and relative subvectors $v_{abs}^{{pt}}$ and
$v_{rel}^{{pt}}$ determined by $I_{pt}$.

We further generalize  the absolute/relative correspondence to the
more general cuts along   submanifolds of arbitrary codimension to
obtain the blow-up/down correspondence. Let $\tilde X$ be a blow-up
of $X$ along a submanifold $Q$ with exceptional divisor $D$. We can
cut $\tilde X$ along $D$ as well as cut $X$ along $Q$. It is
important to observe that the $+$ pairs of these 2 cuts are
essentially the same as the pair $(\tilde X, D)$, in particular,
they have the same relative invariants. Another important fact is
that each invariant of $\tilde X$ in $v_{abs}^{I_{pt}}(\tilde X)$
has a $[pt]$ insertion, and the same is true for $X$. In fact, the
converse is also true. Thus  $\tilde X$  is uniruled if and only if
the absolute vector $v_{abs}^{I_{pt}}(\tilde X)$ is nonzero, and the
same  for $X$.

We now explain why the birational invariance of uniruledness is an
immediate consequence. Suppose the blow-up $\tilde X$ is uniruled,
then $v_{abs}^{{pt}}(\tilde X)$ is nonzero. Hence the relative
vector $v_{rel}^{{pt}}(\tilde X)$ is nonzero by the
absolute/relative correspondence. Apply now the blow-up/down
correspondence to conclude that $v_{abs}^{{pt}}(X)$ is nonzero.
Therefore $X$ is uniruled as well. Similarly we can obtain the
reverse direction.

In the case that the submanifold $D$ is a divisor, another variation
of the absolute/relative correspondence, the so called
sup-admissible correspondence is established in \cite{tLR}. For this
correspondence the relative vector is enlarged to include relative
invariants of $(P_D, D_{\infty})$ with curve classes in the image of
$\iota_*: H_2(D)\to H_2(M)$.

In the case the submanifold $D$ is a uniruled divisor satisfying the
condition of Theorem \ref{main1},  it is further shown in \cite{tLR}
that the sup-admissible correspondence can be restricted  to the
subvector with a relative $[pt]$ insertion and with the curve class
constrained to have symplectic area bounded above by that of a
minimal uniruled class of $D$. The proof  is rather complicated. It
involves the reduction scheme of relative invariants of $\mathbb
P^1-$bundle to invariants of the base in \cite{MP} as well as
\cite{CL} to prove certain vanishing results of relative invariants
of $(P_D, D_{\infty})$.  To prove Theorem \ref{main1},  we also need
a non-vanishing result of relative invariants of $(P_D, D_{\infty})$
to get the divisor/ambient space correspondence.  The outcome of
this correspondence is a nonzero vector of invariants of $X$, each
invariant containing a $[pt]$ insertion. Hence $X$ is uniruled.

\section{Concluding remarks}

Readers can clear sense that the subject of symplectic birational
geometry is only at its beginning.  Constructing uniruled manifolds
by symplectic methods remains to be a challenging problem. We have
already mentioned the issue of singularities: a divisoral
contraction in dimension six already introduces orbifold
singularities. We expect that our program can be carried over to
orbifolds.

 New areas of research include the dichotomy
of higher codimension uniruled submanifolds and transitions. We
ponder whether it makes sense to view transitions as  extended
birational equivalences. Understanding these questions requires new
ideas and technologies.

\subsection{Dichotomy of uniruled submanifold of higher
codimension}  We have a good knowledge of the
       dichotomy of uniruled symplectic divisors at least in dimension six.
       It is natural that we want to expand our understanding to higher codimension uniruled submanifolds.
       There are many reasons to believe that the higher codimension case is very different
       from the divisor case.  In algebraic geometry, this is where we
       encounter other more subtle surgeries such as flip and flop. In the symplectic category, this
       is where Gromov's h-principle is very effective. Therefore, higher codimensional uniruled
       submanifolds should provide a fertile ground for these two completely different theories to interact.

        Corollary \ref{size}  strongly indicates  that
       the {\it size} of a maximal neighborhood   should
       plays an  important role. Such a phenomenon was first observed in  McDuff \cite{Mc3}.
       One could also wonder whether convexity also plays a role
       (compare with \cite{Lai}). It is also desirable to define stable
       and rigid uniruled submanifolds.
       Right now, this is largely an unknown and  exciting territory.


\subsection{Transition}
 Recall
that a transition in the holomorphic category interchanges a
resolution with a smoothing. Symplectically, a smoothing can be
thought of as gluing with a neighborhood of a configuration of
Lagrangian spheres (vanishing cycles). In particular,  a simplest
symplectic transition interchanges a symplectic submanifold with a
Lagrangian sphere. However, a transition is in general not a
birational operation. Thus an important question is to construct
symplectic transitions which are birational cobordism operations.
Such transitions will enhance our ability to `contract' stable
uniruled divisors (submanifolds).


For general transitions it was conjectured by the second author in
\cite{R2} that the quantum cohomology behaves nicely. We could
similarly ask whether uniruledness is preserved under general
transitions.
 Correspondences have been very successful to keep
    track of the total transformation of Gromov-Witten theory
    under birational equivalences.  A natural problem is to construct GW correspondence
    for transitions. This is a new territory for the Gromov-Witten
    theory.

The most famous example is the conifold transition.
    Geometrically, we replace a holomorphic 2-sphere with a
    Lagrangian 3-sphere.
The conifold transition plays an important role in the theory of
Calabi-Yau 3-folds and string theories. In this case one could
partially verify the invariance of uniruledness by \cite{LR}.



     To build a full correspondence
    we may have to enlarge the usual Gromov-Witten theory  to the so
    called open Gromov-Witten theory  to allow  holomorphic curves with
    boundary in Lagrangian manifolds. This is the well-known {\em
    open-closed} duality in physics.
    It also raises an  possibility to extend symplectic
    birational geometry to the {\em open} birational
    geometry. For example, we can define a symplectic manifold to
    be {\em open uniruled} if it contains a nonzero genus
    zero possibly open GW-invariant with a $[pt]$ insertion. Such a notion has already been
    studied in symplectic geometry (cf. \cite{FOOO}).
    Further investigation
    will greatly expand our horizon to understand symplectic
    birational geometry.

\subsection{Final question}
 We feel that  symplectic birational geometry is an interesting
 subject and have raised many questions.
We finish this survey with one more: {\em what kind of structures of
        symplectic manifold are detected by higher genus
       GW invariants}?

\end{document}